\documentclass[12pt, a4paper]{article}
\usepackage{graphicx}
\usepackage{multirow}
\usepackage{color}
\usepackage{tikz}
\usepackage{latexsym}
\usepackage{amssymb}
\usepackage{amsfonts}
\usepackage{amsmath}
\usepackage{indentfirst}
\usepackage{slashbox,colortbl}

\newcommand{\cvd}{\hfill $\blacksquare$\bigskip}
\newtheorem{defn}{Definition}[section]

\newtheorem{prop}{Proposition}[section]

\author{Elena Barcucci\thanks{Dipartimento di Matematica e Informatica ``U. Dini'', Universit\`a degli
		Studi di Firenze, Viale
		G.B. Morgagni 65, 50134 Firenze, Italy. {
			\tt \ elena.barcucci@unifi.it,\quad antonio.bernini@unifi.it,\quad stefano.bilotta@unifi.it,\quad
			renzo.pinzani@unifi.it}}\and Antonio Bernini$^*$ \and Stefano Bilotta$^*$ \and Andrea Lattanzi$^*$ \and Renzo Pinzani$^*$}
		
\title{Non-overlapping Dyck matrices}

\begin{document}
\date{}

\maketitle

\begin{abstract}
We define a set of binary matrices where any two of them can not be placed one on the other in a way such that the corresponding entries coincide. The rows of the matrices are obtained by means of Dyck words. The cardinality of the set of such matrices involves Catalan numbers. 
\end{abstract}


\section{Introduction}

In the last few years several researchers focused their attention to the topic of non-overlapping matrices which can be seen as a generalization of the more famous subject of non-overlapping words. This kind of words is a part of the wide area of Theory of Codes, where, among many interesting questions and applications, the definition of sets of words with particular characteristics is considered. Much care has been paid to matters involving overlaps among the words of a set (or, more precisely, among their prefixes or suffixes). For example, since their introduction, 
various kinds of codes have been considered:  prefix-free, suffix-free and fix-free codes \cite{gilbert,schutz}, comma-free codes \cite{Golomb}, bifix-free codes \cite{nielsen}, non-overlapping codes \cite{black}, cross-bifix-free codes \cite{bajic2,bajic,bilo} (often called non-overlapping codes), variable length non-overlapping codes \cite{bilosolo}.       
More recently, cross-bifix-free codes have been deeply studied 
thanks to their involvement in the development of DNA-based 
storage systems for big amounts of data \cite{DNA1,DNA2}. In 
literature several approaches can be found providing different 
definitions of cross-bifix-free words \cite{succi,bilo,singa} 
and concerning their listing in a Gray code sense \cite{B1,B2}.

The generalization of the notion of non-overlapping strings to two dimensions appears, probably for the first time in a rigorous dissertation, in \cite{tcs} where only square matrices are considered. Here, the notion of prefix or suffix of a string is simply translated into a square matrix anchored in the top-left or bottom-right entry. Nevertheless, in other papers \cite{cripto2D,tcs2} more general matrices are considered and the equivalent of the notion of prefix or suffix is also improved. In particular, a set of non-overlapping matrices is defined over a finite alphabet where the rows of the matrices (which are considered as words) are subjected to some constraints. More precisely, some entries of the matrices are fixed and the rows must avoid certain consecutive patterns of symbols of the alphabet.

\medskip
In the present paper we give a new set of non-overlapping matrices over the binary alphabet. The rows are not restricted words, in the sense that there are no patterns to be avoided. Here, we use the well known Dyck words which allow us to define the matrices of the set in a simpler way with respect to \cite{cripto2D,tcs2}. Moreover, being Dyck words enumerated by the famous Catalan numbers, we are able to find a closed formula for the cardinality of the set, depending only on the dimension of the matrices. A further remark lies in the fact that in each matrix the numbers of $1$'s  and $0$'s are substantially balanced, thanks to the use of Dyck words, differently from the matrices defined in \cite{cripto2D,tcs2}. We conclude by observing that the proofs of the properties of the matrices involve only the examination of prefixes or suffixes of the strings constituting their rows, by using standard techniques typical of non-overlapping words.        

\section{Definitions and background}

The matrices we are dealing with in the present paper are binary matrices and their rows are
generated using particular binary strings (Dyck words). Before describing the construction of the main set (next section), we give a formal definition of non-overlapping matrices, together with some tools useful for the sequel. We also recall some concepts about strings.

\begin{defn}
	Let $u$ and $v$ be two strings over an alphabet
	$\Sigma=\{\sigma_1,\sigma_2\ldots,\sigma_q\}$. They are said \emph{non-overlapping strings} if 
	any non-empty proper prefix of $u$ is different from any non-empty proper suffix of $v$, and vice versa. If $u=v$, then the string $u$ is said \emph{self non-overlapping}.
\end{defn}

If some prefix (suffix) of $u$ matches with some suffix (prefix) of $v$, then $u$ and $v$ are \emph{overlappings strings}. When a factor of length $k$ of consecutive symbols $\sigma_j\in \Sigma$ occurs in a string $u$, we denote it by $\sigma_j^k$. For example, the binary string $u=1000110$ can be indicated as $u=10^31^20$. 

We extend the notion of non-overlapping strings to the case of matrices: 

\begin{defn}\label{def_blocchi}
	Let $A$ and $B$ be two distinct matrices in $M_{m \times n}$, the set of all matrices with $m$ rows and $n$ columns. $A$ and $B$ are said  \emph{non overlapping} if all the following conditions are satisfied by $A$ and $B$: 
	
	\begin{itemize}
		
		\item there do not exist two block partitions\\
				
		$A=\left[\begin{matrix}A_{11} & A_{12} \\ A_{21} & A_{22} \end{matrix}\right]$ \quad  $B=\left[\begin{matrix}B_{11} & B_{12} \\ B_{21} & B_{22} \end{matrix}\right]$\\
		\\
		such that $A_{11}$, $B_{22}$ $\in$ $M_{r \times s}$, with $1\leq r\leq m-1$, $1\leq s\leq n-1$, and either $A_{11}=B_{22}$ nor $A_{12}=B_{21}$, or $A_{21}=B_{12}$, or  $A_{22}=B_{11}$.   \\
		
		\item there do not exist two block partitions \\
		
		$A=\left[\begin{matrix}A_{11} \\ A_{21} \end{matrix}\right]$ \quad  $B=\left[\begin{matrix}B_{11} \\ B_{21} \end{matrix}\right]$\\
		\\
		such that $A_{11}$, $B_{21}$ $\in$ $M_{r \times n}$, with $1\leq r\leq m-1$,
		and either $A_{11}=B_{21}$, or $A_{21}=B_{11}$.\\
		\item   there do not exist two block partitions \\
		
		$A=\left[\begin{matrix}A_{11} & A_{12} \\  \end{matrix}\right]$ \quad  $B=\left[\begin{matrix}B_{11} & B_{12} \\  \end{matrix}\right]$\\
		\\
		such that $A_{11}$, $B_{12}$ $\in$ $M_{m\times s}$, with $1\leq s\leq n-1$, and either $A_{11}=B_{12}$, or  $A_{12}=B_{11}$.\\
	\end{itemize}
	
\end{defn}

\begin{figure}[h!]
	\centering
	\begin{tikzpicture}[scale=0.75]
	\draw (0,0) rectangle (6,4);
	\draw (3,-2) rectangle (9,2);
	\node at (0.5,3.5) {$\mathbf{0}$};
	\node at (1.5,3.5) {$\mathbf{1}$};
	\node at (2.5,3.5) {$\mathbf{0}$};
	\node at (3.5,3.5) {$\mathbf{1}$};
	\node at (4.5,3.5) {$\mathbf{0}$};
	\node at (5.5,3.5) {$\mathbf{1}$};
	\node at (0.5,2.5) {$\mathbf{0}$};
	\node at (1.5,2.5) {$\mathbf{1}$};
	\node at (2.5,2.5) {$\mathbf{1}$};
	\node at (3.5,2.5) {$\mathbf{1}$};
	\node at (4.5,2.5) {$\mathbf{0}$};
	\node at (5.5,2.5) {$\mathbf{1}$};
	\node at (0.5,1.5) {$\mathbf{1}$};
	\node at (1.5,1.5) {$\mathbf{1}$};
	\node at (2.5,1.5) {$\mathbf{1}$};
	\node at (3.5,1.5) {$\mathbf{1}$};
	\node at (4.5,1.5) {$\mathbf{0}$};
	\node at (5.5,1.5) {$\mathbf{0}$};
	\node at (0.5,0.5) {$\mathbf{1}$};
	\node at (1.5,0.5) {$\mathbf{0}$};
	\node at (2.5,0.5) {$\mathbf{0}$};
	\node at (3.5,0.5) {$\mathbf{0}$};
	\node at (4.5,0.5) {$\mathbf{1}$};
	\node at (5.5,0.5) {$\mathbf{1}$};
	\node at (3.5,-1.5) {$\mathbf{0}$};
	\node at (4.5,-1.5) {$\mathbf{1}$};
	\node at (5.5,-1.5) {$\mathbf{0}$};
	\node at (6.5,-1.5) {$\mathbf{1}$};
	\node at (7.5,-1.5) {$\mathbf{0}$};
	\node at (8.5,-1.5) {$\mathbf{1}$};
	\node at (3.5,-0.5) {$\mathbf{1}$};
	\node at (4.5,-0.5) {$\mathbf{1}$};
	\node at (5.5,-0.5) {$\mathbf{1}$};
	\node at (6.5,-0.5) {$\mathbf{0}$};
	\node at (7.5,-0.5) {$\mathbf{0}$};
	\node at (8.5,-0.5) {$\mathbf{0}$};
	\node at (6.5,0.5) {$\mathbf{1}$};
	\node at (7.5,0.5) {$\mathbf{1}$};
	\node at (8.5,0.5) {$\mathbf{0}$};
	\node at (6.5,1.5) {$\mathbf{1}$};
	\node at (7.5,1.5) {$\mathbf{1}$};
	\node at (8.5,1.5) {$\mathbf{1}$};
	\end{tikzpicture}
	\hspace{.5cm}
	\centering
	\begin{tikzpicture}[scale=0.75]
	\draw (0,0) rectangle (6,4);
	\draw (-0.1,-2) rectangle (5.9,2);
	\node at (0.5,3.5) {$\mathbf{0}$};
	\node at (1.5,3.5) {$\mathbf{1}$};
	\node at (2.5,3.5) {$\mathbf{0}$};
	\node at (3.5,3.5) {$\mathbf{1}$};
	\node at (4.5,3.5) {$\mathbf{0}$};
	\node at (5.5,3.5) {$\mathbf{1}$};
	\node at (0.5,2.5) {$\mathbf{0}$};
	\node at (1.5,2.5) {$\mathbf{1}$};
	\node at (2.5,2.5) {$\mathbf{1}$};
	\node at (3.5,2.5) {$\mathbf{1}$};
	\node at (4.5,2.5) {$\mathbf{0}$};
	\node at (5.5,2.5) {$\mathbf{1}$};
	\node at (0.5,1.5) {$\mathbf{1}$};
	\node at (1.5,1.5) {$\mathbf{1}$};
	\node at (2.5,1.5) {$\mathbf{1}$};
	\node at (3.5,1.5) {$\mathbf{1}$};
	\node at (4.5,1.5) {$\mathbf{0}$};
	\node at (5.5,1.5) {$\mathbf{0}$};
	\node at (0.5,0.5) {$\mathbf{1}$};
	\node at (1.5,0.5) {$\mathbf{0}$};
	\node at (2.5,0.5) {$\mathbf{0}$};
	\node at (1.5,-1.5) {$\mathbf{1}$};
	\node at (0.5,-1.5) {$\mathbf{1}$};
	\node at (2.5,-1.5) {$\mathbf{0}$};
	\node at (3.5,-1.5) {$\mathbf{1}$};
	\node at (4.5,-1.5) {$\mathbf{0}$};
	\node at (5.5,-1.5) {$\mathbf{1}$};
	\node at (0.5,-0.5) {$\mathbf{1}$};
	\node at (1.5,-0.5) {$\mathbf{1}$};
	\node at (2.5,-0.5) {$\mathbf{1}$};
	\node at (3.5,-0.5) {$\mathbf{0}$};
	\node at (4.5,-0.5) {$\mathbf{0}$};
	\node at (5.5,-0.5) {$\mathbf{0}$};
	\node at (3.5,0.5) {$\mathbf{1}$};
	\node at (4.5,0.5) {$\mathbf{1}$};
	\node at (5.5,0.5) {$\mathbf{0}$};
	
	\end{tikzpicture}
	
	\vspace{1cm}
	\centering
	\begin{tikzpicture}[scale=0.75]
	\draw (0,0.1) rectangle (6,4.1);
	\draw (4,0) rectangle (10,4);
	\node at (0.5,3.5) {$\mathbf{0}$};
	\node at (1.5,3.5) {$\mathbf{1}$};
	\node at (2.5,3.5) {$\mathbf{0}$};
	\node at (3.5,3.5) {$\mathbf{1}$};
	\node at (4.5,3.5) {$\mathbf{0}$};
	\node at (5.5,3.5) {$\mathbf{1}$};
	
	\node at (0.5,2.5) {$\mathbf{0}$};
	\node at (1.5,2.5) {$\mathbf{1}$};
	\node at (2.5,2.5) {$\mathbf{1}$};
	\node at (3.5,2.5) {$\mathbf{1}$};
	\node at (4.5,2.5) {$\mathbf{0}$};
	\node at (5.5,2.5) {$\mathbf{1}$};
	
	\node at (0.5,1.5) {$\mathbf{1}$};
	\node at (1.5,1.5) {$\mathbf{1}$};
	\node at (2.5,1.5) {$\mathbf{1}$};
	\node at (3.5,1.5) {$\mathbf{1}$};
	\node at (4.5,1.5) {$\mathbf{0}$};
	\node at (5.5,1.5) {$\mathbf{0}$};
	
	\node at (0.5,0.5) {$\mathbf{1}$};
	\node at (1.5,0.5) {$\mathbf{0}$};
	\node at (2.5,0.5) {$\mathbf{0}$};
	\node at (3.5,0.5) {$\mathbf{0}$};
	\node at (4.5,0.5) {$\mathbf{1}$};
	\node at (5.5,0.5) {$\mathbf{1}$};
	

	\node at (6.5,3.5) {$\mathbf{1}$};
	\node at (7.5,3.5) {$\mathbf{1}$};
	\node at (8.5,3.5) {$\mathbf{1}$};
	\node at (9.5,3.5) {$\mathbf{0}$};

	\node at (6.5,2.5) {$\mathbf{1}$};
	\node at (7.5,2.5) {$\mathbf{1}$};
	\node at (8.5,2.5) {$\mathbf{1}$};
	\node at (9.5,2.5) {$\mathbf{0}$};
	
	\node at (6.5,1.5) {$\mathbf{1}$};
	\node at (7.5,1.5) {$\mathbf{0}$};
	\node at (8.5,1.5) {$\mathbf{0}$};
	\node at (9.5,1.5) {$\mathbf{0}$};

	\node at (6.5,0.5) {$\mathbf{0}$};
	\node at (7.5,0.5) {$\mathbf{1}$};
	\node at (8.5,0.5) {$\mathbf{0}$};
	\node at (9.5,0.5) {$\mathbf{1}$};

	\end{tikzpicture}
	\caption{Three pairs of overlapping matrices of $4\times 6$.}
	\label{fig_overlap}
\end{figure}
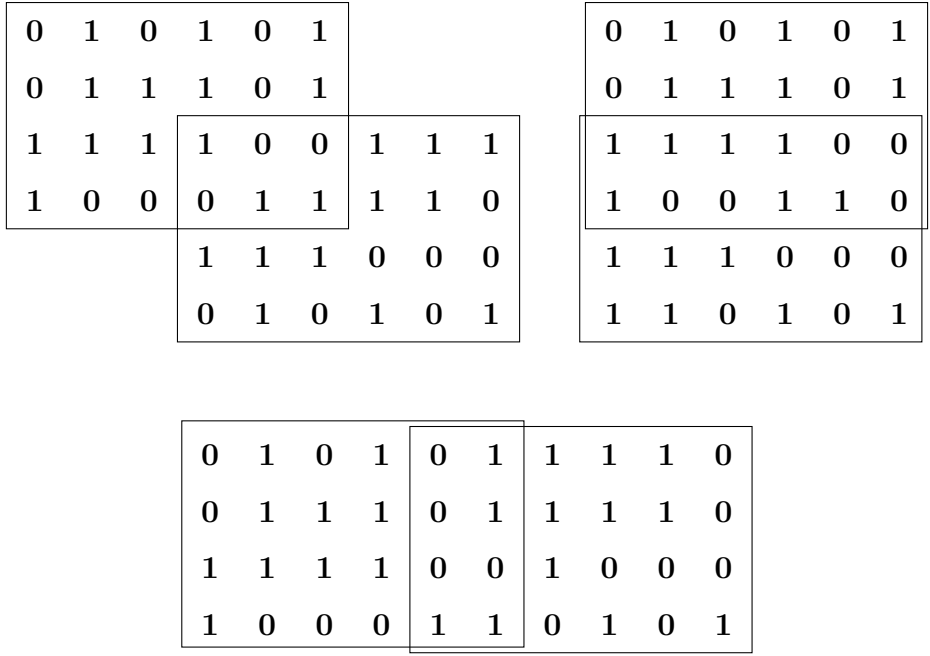

\noindent
In Figure \ref{fig_overlap} three pairs of overlapping matrices contradicting the three above bullets are depicted.

\medskip
\noindent Informally, Definition \ref{def_blocchi} states that two matrices $A$ and $B$ of the same dimension do not overlap if there is no a rigid movement of $B$ on $A$ (note that rotations are not allowed) such that at the end the corresponding entries of $B$ and $A$ coincide. More precisely, after the rigid slipping, a common area (\emph{control window}) is formed. This common area is a rectangular array (the intersection between $A$ and $B$) consisting of a finite number of $1\times 1$ cells of the discrete plane, each containing an entry of $A$ and an entry of $B$. If all the entries of $A$ coincide with the corresponding entries of $B$, then the two matrices are said overlapping matrices, otherwise they are two non-overlapping matrices. Note that the roles of $A$ and $B$ can be interchanged so that we can assume, without loss of generality, that $B$ is the matrix that moves on $A$.
 
Finally, Definition \ref{def_blocchi} can be naturally extended to the case $A=B$, and the matrix is said \emph{self non-overlapping}. 

\begin{defn}
	A set  $L_{m\times n}$ $\subset$ $M_{m\times n}$ is called \emph{non-overlapping}  if each matrix of $L_{m\times n}$ is self non-overlapping and for any two matrices $A,B$ $\in$  $L_{m\times n}$  they are non-overlapping matrices. 
\end{defn}

If $A\in M_{m\times n}$, we indicate, as usual, by $a_{i,j}$ the entry in the $i$-th row and $j$-th column, so that $A=(a_{i,j})$, with $i=1,2,\ldots,m$ and $j=1,2,\ldots,n$. A generic row of $A$ is denoted by $A_i$. Moreover, in the following, if $u$ is a binary string, we denote its length by $|u|$. The number of $1$'s and the number of $0$'s  will be denoted by $|u|_1$ and $|u|_0$, respectively. 

\begin{defn}
	Let $\ell\geqslant 1$. A \emph{Dyck word} $v$ is a binary string of length $2\ell$ such that $|v|_1=|v|_0$ and $|\gamma|_{1}\geq |\gamma|_{0}$ for each prefix $\gamma$ of $v$. The set of all the Dyck words of length $2\ell$ is denoted by $D_{2\ell}$.
\end{defn}

\noindent Note that in a Dyck word it is also $|\delta|_{1}\leqslant |\delta|_{0}$ for any suffix $\delta$ of $v$. For example if $\ell=4$, some Dyck words of length 8 are $11110000, 10101010, 11001100$. 

\begin{defn}
	A Dyck word $v$ is a \emph{type $\alpha$} Dyck word if for each proper prefix $\gamma$, we have $|\gamma|_{1}> |\gamma|_{0}$.
	
\end{defn}

\noindent Clearly, in a \texttt{type} $\alpha$ Dyck word it is also $|\delta|_{1}<|\delta|_{0}$ for any proper suffix $\delta$ of $v$.
Then, a \texttt{type} $\alpha$ Dyck word is a self non-overlapping string \cite{bilo}. Moreover, for each \texttt{type} $\alpha$ Dyck word $v$ of length $2\ell$, it is $v=1u0$ for some Dyck word $u\in D_{2\ell-2}$.   
For example if $\ell=4$, two \texttt{type} $\alpha$ Dyck words of length $8$ are $v=11110000$ or $v=11100100$ where $u=111000$ and $u=110010$, respectively, while the string $10101010$ is not a \texttt{type} $\alpha$ Dyck's word.

\bigskip
We conclude this section by recalling that a Dyck word $u\in D_{2\ell}$ can be read as Dyck path which is a lattice path starting from $(0,0)$, ending at $(2\ell,0)$, never crossing the $x$-axis and using only north-est and south-est steps. A \texttt{type} $\alpha$ Dyck word $v=1u0$, with $|v|=2\ell$ and $u\in D_{2\ell-2}$, is a Dyck path which touches the $x$-axis only at the beginning and at the end. 

\section{A set of non-overlapping binary matrices}

\noindent In this section we are going to define a particular set of non-overlapping binary matrices, the set $L_{m\times n}$. We distinguish two cases depending on the number $n$ of columns: $n$ even and $n$ odd. 

\subsection{The set $L_{m\times n}$ with $n$ even}

\begin{defn}\label{L_m_n_even}
	Let $m\geqslant 2$ and $n\geqslant 4$ with $n$ even. We denote by $L_{m\times n}$ $\subset$ $M_{m\times n}$ the set of matrices  satisfying the following conditions. Let $A=\left(a_{i,j}\right)$
	 an element of $L_{m\times n}$.
	
	\begin{itemize}
		
		\item The first row, $A_{1}$, is a \emph{\texttt{type} $\alpha$} Dyck word of length $n$: $A_1=1u0$ where $u\in D_{n-2}$. Moreover, we require that the first row $A_1$ is the same for all the matrices of $L_{m\times n}$ and it can not appear in any other row.
		
		\item For $i=2,\ldots,m-1$, the row $A_{i}$ is a string of length $n$ chosen among the following different five types of strings: 
		
		\begin{itemize}
			\item \emph{\texttt{type 1}}: a Dyck word of length $n$ ($A_i\in D_n$).
			
			\item \emph{\texttt{type 2}}: a Dyck word of length $n-2$ different from $u$ preceded by $11$ ($A_i=11w$ where $w\in D_{n-2}$ and $w\neq u$).
			
			\item \emph{\texttt{type 3}}: a Dyck word of length $n-2$ different from $u$ followed by 00 ($A_i=w00$ where $w\in D_{n-2}$ and $w\neq u$).
			
			\item \emph{\texttt{type 4}}: a Dyck word of length $n-2$ different from $u$ preceded by $01$ ($A_i=01w$ where $w\in D_{n-2}$ and $w\neq u$).
			
			\item \emph{\texttt{type 5}}: a Dyck word of length $n-2$ followed and preceded by $0$ ($A_i=0w0$ where $w\in D_{n-2}$).
			
			\medskip
			We refer to a row $A_i$ ($i=2,3,\ldots,m-1$) as a \emph{\texttt{type k}} row, with $\emph{\texttt{k}}=1,2,3,4,5$. Note that for \emph{\texttt{type 2, type 3, type 4}} rows it is required $w\neq u$, where $u$ is the Dyck word used to generate $A_1$.
		\end{itemize}
		
		\item The last row $A_{m}$, is a \emph{\texttt{type 1, type 2}} or \emph{\texttt{type 3}} row (\emph{\texttt{type 4}} and \emph{\texttt{type 5}} are not allowed in the last row). 
		
		
	\end{itemize}
	
\end{defn}

 In Figure \ref{7X8} an element of the set $L_{7\times 8}$ is given. The \textit{italic} entries correspond to the ones needed to obtain a \texttt{type 2, type 3, type 4} or \texttt{type 5} row. 
	

	\begin{figure}[h!]
		\centering
		\begin{tikzpicture}[scale=0.75]
		\draw (0,0) rectangle (8,7);
		
		\node at (9.5,6.5) {\texttt{type} $\alpha$ row};
		\node at (0.5,6.5) {$\mathbf{1}$};
		\node at (1.5,6.5) {$\mathbf{1}$};
		\node at (2.5,6.5) {$\mathbf{0}$};
		\node at (3.5,6.5) {$\mathbf{1}$};
		\node at (4.5,6.5) {$\mathbf{0}$};
		\node at (5.5,6.5) {$\mathbf{1}$};
		\node at (6.5,6.5) {$\mathbf{0}$};
		\node at (7.5,6.5) {$\mathbf{0}$};

		\node at (0.5,0.5) {$\mathbf{1}$};
		\node at (1.5,0.5) {$\mathbf{0}$};
		\node at (2.5,0.5) {$\mathbf{1}$};
		\node at (3.5,0.5) {$\mathbf{0}$};
		\node at (4.5,0.5) {$\mathbf{1}$};
		\node at (5.5,0.5) {$\mathbf{0}$};
		\node at (6.5,0.5) {$\mathbf{1}$};
		\node at (7.5,0.5) {$\mathbf{0}$};
		
		\node at (9.5,0.5) {\texttt{type 1} row};
		
		\node at (0.5,3.5) {${\it 0}$};
		\node at (1.5,3.5) {$\it 1$};
		\node at (2.5,3.5) {$\mathbf{1}$};
		\node at (3.5,3.5) {$\mathbf{1}$};
		\node at (4.5,3.5) {$\mathbf{0}$};
		\node at (5.5,3.5) {$\mathbf{0}$};
		\node at (6.5,3.5) {$\mathbf{1}$};
		\node at (7.5,3.5) {$\mathbf{0}$};
		
		\node at (9.5,3.5) {\texttt{type 4} row};

		\node at (0.5,2.5) {$\it {0}$};
		\node at (1.5,2.5) {$\mathbf{1}$};
		\node at (2.5,2.5) {$\mathbf{0}$};
		\node at (3.5,2.5) {$\mathbf{1}$};
		\node at (4.5,2.5) {$\mathbf{0}$};
		\node at (5.5,2.5) {$\mathbf{1}$};
		\node at (6.5,2.5) {$\mathbf{0}$};
		\node at (7.5,2.5) {$\it{0}$};
		
		\node at (9.5,2.5) {\texttt{type 5} row};

		\node at (0.5,1.5) {$\mathbf{1}$};
		\node at (1.5,1.5) {$\mathbf{0}$};
		\node at (2.5,1.5) {$\mathbf{1}$};
		\node at (3.5,1.5) {$\mathbf{0}$};
		\node at (4.5,1.5) {$\mathbf{1}$};
		\node at (5.5,1.5) {$\mathbf{0}$};
		\node at (6.5,1.5) {$\mathbf{1}$};
		\node at (7.5,1.5) {$\mathbf{0}$};
		
		\node at (9.5,1.5) {\texttt{type 1} row};
		
		\node at (0.5,4.5) {$\mathbf{1}$};
		\node at (1.5,4.5) {$\mathbf{0}$};
		\node at (2.5,4.5) {$\mathbf{1}$};
		\node at (3.5,4.5) {$\mathbf{1}$};
		\node at (4.5,4.5) {$\mathbf{0}$};
		\node at (5.5,4.5) {$\mathbf{0}$};
		\node at (6.5,4.5) {$\it{0}$};
		\node at (7.5,4.5) {$\it{0}$};
		
		\node at (9.5,4.5) {\texttt{type 3} row};

		\node at (0.5,5.5) {$\it{1}$};
		\node at (1.5,5.5) {$\it{1}$};
		\node at (2.5,5.5) {$\mathbf{1}$};
		\node at (3.5,5.5) {$\mathbf{1}$};
		\node at (4.5,5.5) {$\mathbf{0}$};
		\node at (5.5,5.5) {$\mathbf{0}$};
		\node at (6.5,5.5) {$\mathbf{1}$};
		\node at (7.5,5.5) {$\mathbf{0}$};
		
		\node at (9.5,5.5) {\texttt{type 2} row};

		\end{tikzpicture}
		\caption{An element of $L_{7\times 8}$.}
		\label{7X8}
	\end{figure}
	

\begin{prop} \label{n_even}
	$L_{m\times n}$ with $n$ even is a non-overlapping set.\\
\end{prop}

\begin{figure}[h!]
	\centering
	\begin{tikzpicture}[scale=0.75]
	\draw (0,0) rectangle (6,4);
	\draw (-0.1,-2) rectangle (5.9,2);
	\node at (0.5,3.5) {$\mathbf{1}$};
	\node at (1.5,3.5) {$\mathbf{1}$};
	\node at (2.5,3.5) {$\mathbf{1}$};
	\node at (3.5,3.5) {$\mathbf{0}$};
	\node at (4.5,3.5) {$\mathbf{0}$};
	\node at (5.5,3.5) {$\mathbf{0}$};
	
	\node at (0.5,1.5) {$\mathbf{1}$};
	\node at (1.5,1.5) {$\mathbf{1}$};
	\node at (2.5,1.5) {$\mathbf{1}$};
	\node at (3.5,1.5) {$\mathbf{0}$};
	\node at (4.5,1.5) {$\mathbf{0}$};
	\node at (5.5,1.5) {$\mathbf{0}$};
	\end{tikzpicture}
	\caption{Vertical overlap}
	\label{vertical overlap}
\end{figure}
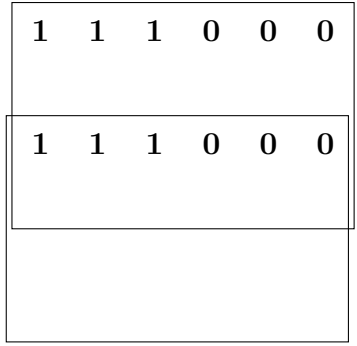

\noindent \emph{Proof.} Let $A$ and $B$ two matrices in $L_{m\times n}$ (possibly the same). Since we have to check if $A$ and $B$ overlap, we immediately notice that the first row of $A$ or $B$ is always involved in the control window. We consider two different kinds of overlaps: \emph{vertical overlap} when the $i$-column of $A$ is superimposed to the $i$-column of $B$ and \emph{skew overlap} in all the other cases.\\
If a vertical overlap occurs, then $A$ and $B$ can not overlap since, looking at the control window (Figure \ref{vertical overlap}), the row $A_1$ ($B_1=A_1$) is over a row $B_j$ ($A_j$) and, by Definition \ref{L_m_n_even}, $A_1\neq B_j$ ($B_1\neq A_j$), for $j=2,3,\ldots,m$.

\medskip
If $A$ and $B$ skew overlap, then two cases must be considered: the control window contains a prefix of $A_1$ (or $B_1=A_1)$ or the control windows contains a suffix of $A_1$ or $B_1$.

\begin{itemize}
	\item Let $\gamma$ be a proper prefix of $A_1$ which is involved in the control window: it is superimposed with a proper suffix $\delta$ of $B_j$, $j=1,2,\ldots,m$. Depending on the type of the row $B_j$, we have different cases:
\begin{itemize}
	
	\item $B_j$ is a \texttt{type 1} row: since $|\gamma|_1>|\gamma|_0$ and 
	$|\delta|_0\geqslant |\delta|_1$, the prefix $\gamma$ and the suffix $\delta$ do not match, then rows $A_1$ and $B_j$, and consequently the matrices $A$ and $B$, do not overlap.	
	
	\item $B_j$ is a \texttt{type 2} or a \texttt{type 4} row: if $|\delta|\leqslant n-2$, then the 
	same argument of the above case can be used to conclude that $A$ and $B$ do 
	not overlap, while, if $|\delta|=n-1$, then we observe that
	$\delta=1 b_{j,3}b_{j,4}\ldots b_{j,n}$ and $\gamma=1u$ where $u\in D_{n-2}$, so that $\gamma$ and $\delta$ do not match since, by Definition \ref{L_m_n_even}, it is $u\neq b_{j,3}b_{j,4}\ldots b_{j,n}$. Therefore, $A$ and $B$ do not overlap.	
	

%
	\item $B_j$ is a \texttt{type 3} or a \texttt{type 5} row: since $|\gamma|_1>|\gamma|_0$ and $|\delta|_0 > |\delta|_1$, matrices $A$ and $B$ can not overlap.

%
	
\end{itemize}

Notice that in the above bullet, the roles of $A$ and $B$ can be interchanged, leading to the same conclusion.

\item Let $\delta$ be a proper suffix of $A_1$ which is involved in the control window: it is superimposed with a proper prefix $\gamma$ of $B_j$, $j=1,2,\ldots,m$.
Depending on the type of the row $B_j$, we have different cases:


\begin{itemize}
	
	\item $B_j$ is a \texttt{type 1} row: since $|\delta|_0>|\delta|_1$ and $|\gamma|_1\geqslant |\gamma|_0$, the prefix $\gamma$ and the suffix $\delta$ do not match, then rows $A_1$ and $B_j$, and consequently the matrices $A$ and $B$, do not overlap.
	
	\item $B_j$ is a \texttt{type 2} row: since $|\delta|_0>|\delta|_1$ and $|\gamma|_1> |\gamma|_0$, similarly to the above case, matrices $A$ and $B$ can not overlap.

	
	\item $B_j$ is a \texttt{type 3} row: if $|\gamma|\leqslant n-2$, no match is possible (see before) between $\gamma$ and $\delta$ (so that $A$ and $B$ do not overlap), while, if $|\gamma|= n-1$, then we observe that $\gamma=b_{j,1}b_{j,2}\ldots b_{j,n-2}0$ and $\delta=u0$ where $u\in D_{n-2}$, so that $\gamma$ and $\delta$ do not match since, by Definition \ref{L_m_n_even}, it is $u\neq b_{j,1}b_{j,2}\ldots b_{j,n-2}$. Therefore $A$ and $B$ can not overlap.

	
	\item $B_j$ is a \texttt{type 4} row: in this case, by Definition \ref{L_m_n_even}, it is $j\neq m$. If $|\gamma|=1$, then $\gamma=b_{1,j}=0$ and $\delta=a_{1,n}=0$ (each row of a matrix in the set ends with 0) form a match. Nevertheless, there exists an index $i$ (with $2\leqslant i\leqslant m-1$) such that $b_{1,m}=1$ is over the entry $a_{i,n}$ which is 0. Therefore, no overlapping between $A$ and $B$ is possible (Figure \ref{FIG_X}). If $|\gamma|>1$, then $|\gamma|_1\geqslant |\gamma|_0$ and $|\delta|_0>|\delta|_1$ which does not allow $A$ and $B$ to overlap. 


	\item $B_j$ is a \texttt{type 5} row: also in this case $j\neq m$. If $|\gamma|=1$, the matrices $A$ and $B$ are arranged as in the above case and they do not overlap (Figure \ref{FIG_X}).
	If $|\gamma| >1$, then $\gamma=0w$ where $w$ is a proper prefix of a Dyck word, so that $|\gamma|_1\geqslant |\gamma|_0$ (if $|w|_1 >|w|_0$) or $|\gamma|_1+1=|\gamma|_0$ (if $|w|_1=|w|_0$). In the first case $\gamma$ and $\delta$ do not match since it is $|\delta|_0>|\delta|_1$. In the latter, the string $w$ would be a suffix of $A_1$. This is not possible since $A_1$ is a \texttt{type} $\alpha$ row. We conclude that also in this case $A$ and $B$ do not overlap.

\bigskip
	
	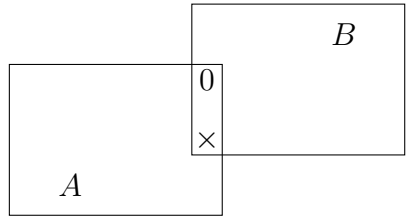
\begin{figure}[h!]
		\centering
		\begin{tikzpicture}[scale=0.4]
		\draw (0,0) rectangle (7,5);
		\draw (6,2) rectangle (13,7);
		\node at (6.5,4.5) {0};
		\node at (6.5,2.5) {$\times$};
		\node at  (2,1) {$A$};
		\node at  (11,6) {$B$};

		\end{tikzpicture}
		\caption{The prefix of a type 4 or type 5 row in $B$ composed by the only first 0  is superimposed on the last 0 of the first row in $A$. The $\times$ represents a mismatch.}
	\label{FIG_X}	
\end{figure}

\end{itemize}

\end{itemize}

Notice that from the above proof we can deduce that $A_1$ and any other row of \texttt{type 1, 2} or \texttt{3} are non-overlapping strings. 

\cvd

\subsection{The set $L_{m\times n}$ with $n$ odd}

\begin{defn}\label{L_m_n_odd}
	Let $m\geqslant 2$ and $n\geqslant 5$ with $n$ odd. We denote by $L_{m\times n}$ $\subset$ $M_{m\times n}$ the set of matrices satisfying the following conditions. Let $A=\left(a_{i,j}\right)$ an element of $L_{m\times n}$.
	
	\begin{itemize}
		
		\item The first row, $A_{1}$, is a \emph{\texttt{type}} $\alpha$ Dyck word of length $n-1$ preceded by 1: $A_1=11u0$ where $u\in D_{n-3}$.  Moreover, we require that the first row is the same for all the matrices of $L_{m\times n}$ and it can not appear in any other row. 
		
		\item For $i=2,\ldots,m-1$, the row $A_{i}$ is a string of length $n$ chosen among the following different three types of strings: 
		
		\begin{itemize}
						
			\item \emph{\texttt{type 6}}: a Dyck word of length $n-1$ preceded by 1 ($A_i=1w$ where $w\in  D_{n-1}$ and $w\neq u$).
			
			\item \emph{\texttt{type 7}}: a Dyck word of length $n-1$ followed by 0 ($A_i=w0$ where $w\in  D_{n-1}$ and $w\neq u$).
			
			\item \emph{\texttt{type 8}}: a Dyck word of length $n-1$ preceded by 0 ($A_i=0w$ where $w\in  D_{n-1}$).
						
			\medskip
			We refer to a row $A_i$ ($i=2,3,\ldots,m-1$) as a \emph{\texttt{type k}} row, with $\emph{\texttt{k}}=6,7,8$. Note that for \emph{\texttt{type 6, type 7}} rows it is required $w\neq u$, where $u$ is the Dyck word used to generate $A_1$.
		\end{itemize}
		
		\item The last row $A_{m}$ is a \emph{\texttt{type 6}} or \emph{\texttt{type 7}} row (\emph{\texttt{type 8}} is not allowed in the last row).

	\end{itemize}
	
\end{defn}

 In Figure \ref{6X7} an element of the set $L_{6\times 7}$ is given. The \textit{italic} entries correspond to the ones needed to obtain a \texttt{type 6, type 7} or \texttt{type 8} row.

	\begin{figure}[h!]
		\centering
		\begin{tikzpicture}[scale=0.75]
		\draw (0,0) rectangle (7,6);
		
		\node at (0.5,5.5) {$\it{1}$};
		\node at (1.5,5.5) {$\mathbf{1}$};
		\node at (2.5,5.5) {$\mathbf{1}$};
		\node at (3.5,5.5) {$\mathbf{1}$};
		\node at (4.5,5.5) {$\mathbf{0}$};
		\node at (5.5,5.5) {$\mathbf{0}$};
		\node at (6.5,5.5) {$\mathbf{0}$};
		
		\node at (8.1,5.5) {\texttt{type} $\alpha$};

		\node at (0.5,0.5)  {$\mathbf{1}$};
		\node at (1.5,0.5) {$\mathbf{1}$};
		\node at (2.5,0.5) {$\mathbf{0}$};
		\node at (3.5,0.5) {$\mathbf{0}$};
		\node at (4.5,0.5) {$\mathbf{1}$};
		\node at (5.5,0.5) {$\mathbf{0}$};
		\node at (6.5,0.5) {$\it{0}$};
		
		\node at (8.1,0.5) {type 7};
		
		\node at (0.5,3.5) {$\mathbf{1}$};
		\node at (1.5,3.5) {$\mathbf{0}$};
		\node at (2.5,3.5) {$\mathbf{1}$};
		\node at (3.5,3.5) {$\mathbf{0}$};
		\node at (4.5,3.5) {$\mathbf{1}$};
		\node at (5.5,3.5) {$\mathbf{0}$};
		\node at (6.5,3.5) {$\it{0}$};
		
		\node at (8.1,3.5) {type 7};
		
		\node at (0.5,2.5) {$\it{0}$};
		\node at (1.5,2.5) {$\mathbf{1}$};
		\node at (2.5,2.5) {$\mathbf{1}$};
		\node at (3.5,2.5) {$\mathbf{1}$};
		\node at (4.5,2.5) {$\mathbf{0}$};
		\node at (5.5,2.5) {$\mathbf{0}$};
		\node at (6.5,2.5) {$\mathbf{0}$};
		
		\node at (8.1,2.5) {type 8};
		
		\node at (0.5,1.5) {$\it{1}$};
		\node at (1.5,1.5) {$\mathbf{1}$};
		\node at (2.5,1.5) {$\mathbf{1}$};
		\node at (3.5,1.5) {$\mathbf{0}$};
		\node at (4.5,1.5) {$\mathbf{0}$};
		\node at (5.5,1.5) {$\mathbf{1}$};
		\node at (6.5,1.5) {$\mathbf{0}$};
		
		\node at (8.1,1.5) {type 6};
		
		\node at (0.5,4.5) {$\it{1}$};
		\node at (1.5,4.5) {$\mathbf{1}$};
		\node at (2.5,4.5) {$\mathbf{0}$};
		\node at (3.5,4.5) {$\mathbf{1}$};
		\node at (4.5,4.5) {$\mathbf{0}$};
		\node at (5.5,4.5) {$\mathbf{1}$};
		\node at (6.5,4.5) {$\mathbf{0}$};
		
		\node at (8.1,4.5) {type 6};

		\end{tikzpicture}
		\caption{An element of $L_{6\times 7}$.}
		\label{6X7}
	\end{figure}

\begin{prop}\label{n_odd}
	$L_{m\times n}$ with $n$ odd is a non-overlapping set.\\
\end{prop}

\noindent \emph{Proof.} Let $A$ and $B$ two matrices in $L_{m\times n}$ (possibly the same). The argument is very similar to the one used in the proof of the Proposition \ref{n_even}. 
 Since we have to check if $A$ and $B$ overlap, we immediately notice that the first row of $A$ or $B$ is always involved in the control window.

 We distinguish two different kinds of overlaps: \emph{vertical overlap} and \emph{skew overlap}.
If a vertical overlap occurs, then $A$ and $B$ can not overlap since, looking at the control window, the row $A_1$ ($B_1$) is over a row $B_j$ ($A_j$) and, by Definition \ref{L_m_n_odd}, $A_1\neq B_j$ ($B_1\neq A_j$), for $j=2,3,\ldots,m$.  
If $A$ and $B$ skew overlap, then two cases must be considered: the control window contains a prefix of $A_1$ or $B_1$ or the control windows contains a suffix of $A_1$ or $B_1$.

\begin{itemize}
	\item Let $\gamma$ be a proper prefix of $A_1$ which is involved in the control window: it is superimposed with a proper suffix $\delta$ of $B_j$, $j=1,2,\ldots,m$. Depending on the type of the row $B_j$, we have different cases:
	\begin{itemize}
		
		\item $B_j$ is a \texttt{type 6} or a \texttt{type 8} row: since $|\gamma|_1>|\gamma|_0$ and 
		$|\delta|_0\geqslant |\delta|_1$, the prefix $\gamma$ and the suffix $\delta$ do not match, then rows $A_1$ and $B_j$, and consequently the matrices $A$ and $B$, do not overlap.	
		
		\item $B_j$ is a \texttt{type 7} row: since $|\gamma|_1>|\gamma|_0$ and 
		$|\delta|_0 > |\delta|_1$, the prefix $\gamma$ and the suffix $\delta$ do not match, then rows $A_1$ and $B_j$, and consequently the matrices $A$ and $B$, do not overlap.	
		
	\end{itemize}
	
	Notice that in the above bullet, the roles of $A$ and $B$ can be interchanged, leading to the same conclusion.

	\item Let $\delta$ be a proper suffix of $A_1$ which is involved in the control window: it is superimposed with a proper prefix $\gamma$ of $B_j$, $j=1,2,\ldots,m$.
	Depending on the type of the row $B_j$, we have different cases:
	
	\begin{itemize}
		
		\item $B_j$ is a \texttt{type 6} row: since $|\delta|_0\geqslant |\delta|_1$ ($|\delta|_0> |\delta|_1$ if $|\delta|\leqslant n-2)$ and $|\gamma|_1> |\gamma|_0$, the prefix $\gamma$ and the suffix $\delta$ do not match, then rows $A_1$ and $B_j$, and consequently the matrices $A$ and $B$, do not overlap.
		
		\item $B_j$ is a \texttt{type 7} row: if $|\delta|\leqslant n-2$, since $|\delta|_0 > |\delta|_1$ and $|\gamma|_1 \geqslant |\gamma|_0$, the prefix $\gamma$ and the suffix $\delta$ do not match, then rows $A_1$ and $B_j$, and consequently the matrices $A$ and $B$, do not overlap.\\
		\noindent If $|\delta|=n-1$, then $\delta =u$, where $u \in D_{n-1}$ and $\gamma = b_{j,1}b_{j,2}\ldots b_{j,n-1} =w$ where $w \in D_{n-1}$ so that $\gamma$ and $\delta$ do not match since, by Definition \ref{n_odd}, it is $u\neq w$. Therefore $A$ and $B$ can not overlap.


		\item $B_j$ is a \texttt{type 8} row: in this case $j\neq m$. 
		If $|\gamma|=1$, then $\gamma=b_{1,j}=0$ and $\delta=a_{1,n}=0$ (each row of a matrix in the set ends with 0) form a match. Nevertheless, there exists an index $i$ (with $2\leqslant i\leqslant m-1$) such that $b_{1,m}=1$ is over the entry $a_{i,n}$ which is 0. Therefore, no overlapping between $A$ and $B$ is possible (Figure \ref{FIG_X}). 	If $|\gamma| >1$, then $\gamma=0w$ where $w$ is a proper prefix of a Dyck word, so that $|\gamma|_1\geqslant |\gamma|_0$ (if $|w|_1 >|w|_0$) or $|\gamma|_1+1=|\gamma|_0$ (if $|w|_1=|w|_0$). In the first case $\gamma$ and $\delta$ do not match since it is $|\delta|_0>|\delta|_1$. In the latter, the string $w$ would be a suffix of $A_1$. This is not possible since $A_1$ is a \texttt{type} $\alpha$ row. We conclude that also in this case $A$ and $B$ do not overlap. 
		
	\end{itemize}

\end{itemize}

Notice that, similarly to Proposition \ref{n_even}, from the above proof we can deduce that $A_1$ and any other row of \texttt{type 6} or \texttt{7} are non-overlapping strings. 

\cvd

	\section{The enumeration of  $L_{m\times n}$}
	
	In this section we are going to enumerate the set  $L_{m\times n}$. The first row $A_1$ of all the matrices in $L_{m\times n}$ is the same, so that it is easy to realize that its cardinality depends on the number of the remaining $m-1$ rows satisfying the constraints of Definition \ref{L_m_n_even} or \ref{L_m_n_odd}. We recall that the Dyck words of length $2n$ are enumerated by the $n$-th Catalan number $C_n$ given by
		
		$$C_n=\dfrac{1}{n+1} \binom{2n}{n}\ .$$
		
%
	Following the constraints of Definition \ref{L_m_n_even}, in the case  $n$ even and reminding that $A_1=1u0$ with $u\in D_{n-2}$, we observe that if the row $A_i$, $i=2,3,\ldots,m-1$, is a \texttt{type 2,3,4} row, then it can be constructed by using a Dyck word $w\in D_{n-2}$ and $w\neq u$ and adding the prefixes $11$ or $01$ or the suffix $00$, so that there are $C_{\frac{n-2}{2}}-1$ possibilities for each type of row. If $A_i$ is a \texttt{type 5} row, then it can be constructed in $C_{\frac{n-2}{2}}$ different ways by using whichever Dyck word $w\in D_{n-2}$ and adding the prefix and the suffix $0$ ($A_i=0w0$). Finally, if $A_i$ is a \texttt{type 1} row, then $A_i\in D_n$ and $A_i\neq A_1$. The foregoing is summarized in the first term of the following formula where the exponent $m-2$ tracks the number of rows from the second to the second to last one. The second term counts the possibilities for the last row which is a \texttt{type 1, 2} or \texttt{type 3} row.

	$$|L_{m \times n}|=\left(\underbrace{ C_{\frac{n}{2}}-1}_{\texttt{type} \ 1} + \underbrace{3C_{\frac{n-2}{2}}-3}_{\texttt{type} \ \texttt{2},\ \texttt 3,\ \texttt 4} +\underbrace{C_{\frac{n-2}{2}}}_{\texttt{type} \ \texttt 5} \right)^{(m-2)} \left(  \underbrace{C_{\frac{n}{2}}-1}_{\texttt{type} \ \texttt 1} + \underbrace{2C_{\frac{n-2}{2}}-2}_{\texttt{type}\ \texttt 2,\ \texttt 3}   \right)$$


	 With a similar argument, following Definition \ref{L_m_n_odd} in the case $n$ odd, we obtain:

	$$|L_{m \times n}|=\left(\underbrace{2C_{\frac{n-1}{2}}-2}_
	{\texttt{type}  \ \texttt 6,\ \texttt 7  }+ \underbrace{C_{\frac{n-1}{2}}}_{\texttt{type} \ \texttt 8} \right)^{(m-2)}
	\left(\underbrace{2C_{\frac{n-2}{2}}-2}_{\texttt{type}  \ \texttt 6,\ \texttt 7} \right)$$
	
	\bigskip
	\noindent Table \ref{values} shows some values of the cardinality $L_{m \times n}$ for some estimated values of $m$ and $n$ up to 10. A careful analysis of them, reveals that if $n>m$, then the cardinality of $L_{m\times n}$ is greater then the cardinality of $L_{n\times m}$. So, fixed $n$ and $m$, in the case $n>m$, it is convenient to define the matrices by columns (according to Definition \ref{L_m_n_even} and \ref{L_m_n_odd}) in order to obtain the set with the highest cardinality. In other words, if the columns are longer than the rows, then Dyck words should be used to construct the columns of the matrices of the set. 
	\begin{table}
	\begin{center}	
	$$\begin{tabular}{|r|c|c|c|c|c|c|c|}
	\hline
	$ m \diagdown n $  &4 & 5 & 6   & 7 &8 &9 &10 \\ \hline             
	2 &1 & 2 & 6& 8&  21 &  26  & 67\\ \hline
	3& 2& 4 &54 & 104 &  630  &  1040  & $6.2\cdot 10^{3}$\\ \hline
	
	4 & 4 & 8 & 486& 1352 &  $1.9\cdot 10^{4}$   & $4.1\cdot 10^{4}$   & $5.9\cdot 10^{5}$\\ \hline
	
	5&8 &  16& 4374  &  $1.7\cdot 10^{4}$ &  $5.7\cdot 10^{5}$   &  $1.6\cdot 10^{6}$  &  $5.5\cdot 10^{7}$\\ \hline
	6& 16 &   32  & $3.9\cdot 10^{4}$  & $2.2\cdot 10^{5}$   &  $1.7\cdot 10^{7}$   &  $6.6\cdot 10^{7}$  & $5.2\cdot 10^{9}$\\ \hline
	7&32 &   64  &$3.5\cdot 10^{5}$ &  $3.0\cdot 10^{6}$  &  $5.1\cdot 10^{8}$   & $2.7\cdot 10^{9}$   &  $4.9\cdot 10^{11}$\\ \hline
	8&64 &   128   & $3.1\cdot 10^{6}$  &  $3.8\cdot 10^{7}$  & $1.5\cdot 10^{10}$ &  $1.1\cdot 10^{11}$  & $4.6\cdot 10^{13}$\\ \hline
	9& 128 &  256   &  $2.8\cdot 10^{7}$  &  $5.0\cdot 10^{8}$  &   $4.6\cdot 10^{11}$ &  $4.2\cdot 10^{12}$  & $4.3\cdot 10^{15}$ \\ \hline
	10& 256 &512 &  $2.6\cdot 10^{8}$  &  $6.5\cdot 10^{9}$  &  $1.4\cdot 10^{13}$  & $1.7\cdot 10^{14}$   & $4.1\cdot 10^{17}$\\ \hline
	\end{tabular}$$\\
	\caption{Some values of the cardinality of $L_{m\times n}$}
	\label{values}
\end{center}	
\end{table}
	
\section{On the non-expandability of $L_{m\times n}$}

As well as in the theory of  
cross-fix-free codes, an interesting study concerns their 
non-expandability. We recall the following definition, adapted to matrices.

\begin{defn}
	
	Let $N$ be a set of non-overlapping matrices and let $U$ the set of all self-non-overlapping matrices. $N$ is  \emph{non-expandable} if for each $A$ $\in$ $U \setminus N$  there exists at least one matrix $B \in N$ such that $A$ and $B$ are not non-overlapping matrices.

\end{defn}

\noindent On the contrary, the set $N$ is \emph{expandable} if 
there exists a matrix $A \in U \setminus N$ such that $A$ and $B$ are non-overlapping matrices, for any $B\in N$.\\

We will prove that $L_{m\times n}$ is not a non-expandable set of matrices. The reason why it is not, maybe lies on the fact that, for such kinds of sets, some constraints on the columns of the matrices are needed, as it seems to come out from \cite{NENO}.
The columns of the matrices of $L_{m\times n}$ are not subject to obvious constraints. We recall that $L_{m\times n}$ is constructed by using Dyck words which form the rows of the matrices (adding some suitable prefixes or suffixes). This does not induces explicit restrictions over the columns.

\medskip
Fixed a set $L_{m\times n}$, denote $A_1$ the first row of 
each $A\in L_{m\times n}$. Suppose that there exist a binary 
strings $X$ with the following properties:
\begin{enumerate}
	\item $X$ is not a \texttt{type 1, 2, 3, 4} or \texttt{type 5} string in the $n$ even case (\texttt{type 6, 7} or \texttt{type 8} in the $n$ odd case);

	\item $X$ and $A_1$ are non-overlapping strings.
\end{enumerate}

\noindent
We observe that, from the second property, $X$ begins with $1$ and ends with $0$.
Now, we define the matrix $Z\in M_{n\times n}$ as follows:	
\begin{itemize}
	\item $Z_{1}=A_1$; 
	\item $Z_{i}=X$\quad  for $i=2,\ldots , m$.
\end{itemize}

\noindent It is not difficult to realize that $Z\notin L_{m\times n}$: its rows from the second to the last one do not match the requirements of Definitions \ref{L_m_n_even} or \ref{L_m_n_odd}. Moreover, $Z$ is a self non-overlapping matrix
since any prefix (suffix) of $Z_1$ (involved in any vertical or skew overlap) does not match any suffix (prefix) of $Z_i=X$ 
for the second property of the string $X$.

We conclude with the following

\begin{prop}\label{expanda}
	
	$Z \cup L_{m\times n}$ is a non-overlapping set. 
	
\end{prop}

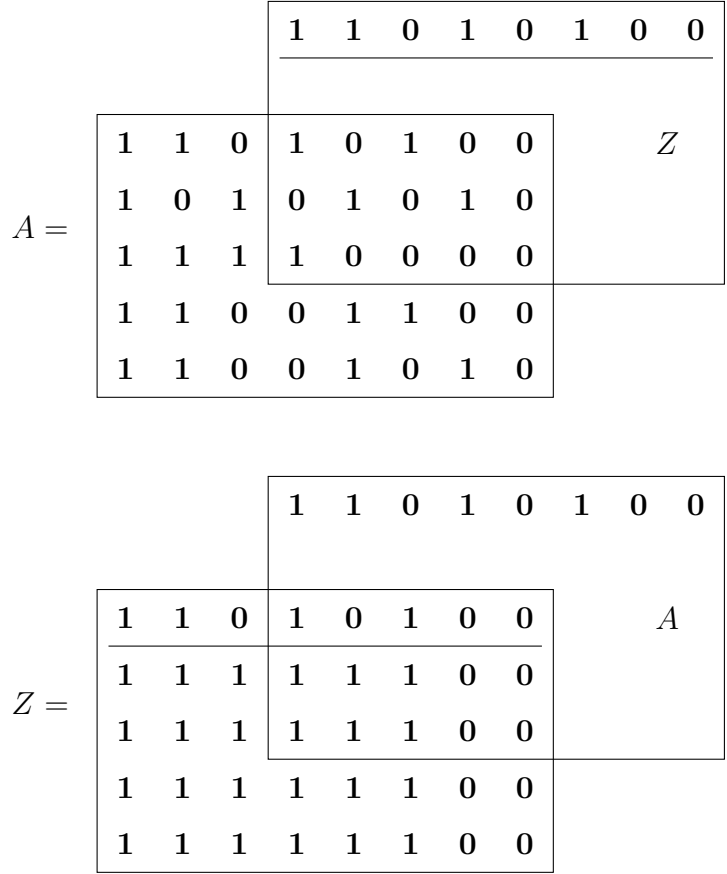
\begin{figure}
	\centering
	\begin{tikzpicture}[scale=0.75]
	
	\draw (0,0) rectangle (8,5);
	\draw (3,2) rectangle (11,7);

	\node at (-1,3) {$A=$};
	\node at (10,4.5) {$Z$};

	\node at (3.5,6.5) {$\mathbf{1}$};
	\node at (4.5,6.5) {$\mathbf{1}$};
	\node at (5.5,6.5) {$\mathbf{0}$};
	\node at (6.5,6.5) {$\mathbf{1}$};
	\node at (7.5,6.5) {$\mathbf{0}$};
	\node at (8.5,6.5) {$\mathbf{1}$};
	\node at (9.5,6.5) {$\mathbf{0}$};
	\node at (10.5,6.5) {$\mathbf{0}$};

	\node at (0.5,4.5) {$\mathbf{1}$};
	\node at (1.5,4.5) {$\mathbf{1}$};
	\node at (2.5,4.5) {$\mathbf{0}$};
	\node at (3.5,4.5) {$\mathbf{1}$};
	\node at (4.5,4.5) {$\mathbf{0}$};
	\node at (5.5,4.5) {$\mathbf{1}$};
	\node at (6.5,4.5) {$\mathbf{0}$};
	\node at (7.5,4.5) {$\mathbf{0}$};
	
	\node at (0.5,3.5) {$\mathbf{1}$};
	\node at (1.5,3.5) {$\mathbf{0}$};
	\node at (2.5,3.5) {$\mathbf{1}$};
	\node at (3.5,3.5) {$\mathbf{0}$};
	\node at (4.5,3.5) {$\mathbf{1}$};
	\node at (5.5,3.5) {$\mathbf{0}$};
	\node at (6.5,3.5) {$\mathbf{1}$};
	\node at (7.5,3.5) {$\mathbf{0}$};
	
	\node at (0.5,2.5) {$\mathbf{1}$};
	\node at (1.5,2.5) {$\mathbf{1}$};
	\node at (2.5,2.5) {$\mathbf{1}$};
	\node at (3.5,2.5) {$\mathbf{1}$};
	\node at (4.5,2.5) {$\mathbf{0}$};
	\node at (5.5,2.5) {$\mathbf{0}$};
	\node at (6.5,2.5) {$\mathbf{0}$};
	\node at (7.5,2.5) {$\mathbf{0}$};
	
	\node at (0.5,1.5) {$\mathbf{1}$};
	\node at (1.5,1.5) {$\mathbf{1}$};
	\node at (2.5,1.5) {$\mathbf{0}$};
	\node at (3.5,1.5) {$\mathbf{0}$};
	\node at (4.5,1.5) {$\mathbf{1}$};
	\node at (5.5,1.5) {$\mathbf{1}$};
	\node at (6.5,1.5) {$\mathbf{0}$};
	\node at (7.5,1.5) {$\mathbf{0}$};
	
	\draw (3.2,6) -- (10.8,6);
	
	\node at (0.5,0.5) {$\mathbf{1}$};
	\node at (1.5,0.5) {$\mathbf{1}$};
	\node at (2.5,0.5) {$\mathbf{0}$};
	\node at (3.5,0.5) {$\mathbf{0}$};
	\node at (4.5,0.5) {$\mathbf{1}$};
	\node at (5.5,0.5) {$\mathbf{0}$};
	\node at (6.5,0.5) {$\mathbf{1}$};
	\node at (7.5,0.5) {$\mathbf{0}$};
	\end{tikzpicture}
	
	\vspace{1cm}
	
	\centering
	\begin{tikzpicture}[scale=0.75]
	
	\draw (0,0) rectangle (8,5);
	\draw (3,2) rectangle (11,7);

	\node at (-1,3) {$Z=$};
	\node at (10,4.5) {$A$};

	\node at (3.5,6.5) {$\mathbf{1}$};
	\node at (4.5,6.5) {$\mathbf{1}$};
	\node at (5.5,6.5) {$\mathbf{0}$};
	\node at (6.5,6.5) {$\mathbf{1}$};
	\node at (7.5,6.5) {$\mathbf{0}$};
	\node at (8.5,6.5) {$\mathbf{1}$};
	\node at (9.5,6.5) {$\mathbf{0}$};
	\node at (10.5,6.5) {$\mathbf{0}$};

	\node at (0.5,4.5) {$\mathbf{1}$};
	\node at (1.5,4.5) {$\mathbf{1}$};
	\node at (2.5,4.5) {$\mathbf{0}$};
	\node at (3.5,4.5) {$\mathbf{1}$};
	\node at (4.5,4.5) {$\mathbf{0}$};
	\node at (5.5,4.5) {$\mathbf{1}$};
	\node at (6.5,4.5) {$\mathbf{0}$};
	\node at (7.5,4.5) {$\mathbf{0}$};
	
	\node at (0.5,3.5) {$\mathbf{1}$};
	\node at (1.5,3.5) {$\mathbf{1}$};
	\node at (2.5,3.5) {$\mathbf{1}$};
	\node at (3.5,3.5) {$\mathbf{1}$};
	\node at (4.5,3.5) {$\mathbf{1}$};
	\node at (5.5,3.5) {$\mathbf{1}$};
	\node at (6.5,3.5) {$\mathbf{0}$};
	\node at (7.5,3.5) {$\mathbf{0}$};
	
	\node at (0.5,2.5) {$\mathbf{1}$};
	\node at (1.5,2.5) {$\mathbf{1}$};
	\node at (2.5,2.5) {$\mathbf{1}$};
	\node at (3.5,2.5) {$\mathbf{1}$};
	\node at (4.5,2.5) {$\mathbf{1}$};
	\node at (5.5,2.5) {$\mathbf{1}$};
	\node at (6.5,2.5) {$\mathbf{0}$};
	\node at (7.5,2.5) {$\mathbf{0}$};
	
	\node at (0.5,1.5) {$\mathbf{1}$};
	\node at (1.5,1.5) {$\mathbf{1}$};
	\node at (2.5,1.5) {$\mathbf{1}$};
	\node at (3.5,1.5) {$\mathbf{1}$};
	\node at (4.5,1.5) {$\mathbf{1}$};
	\node at (5.5,1.5) {$\mathbf{1}$};
	\node at (6.5,1.5) {$\mathbf{0}$};
	\node at (7.5,1.5) {$\mathbf{0}$};
	
	\draw (0.2,4) -- (7.8,4);
	
	\node at (0.5,0.5) {$\mathbf{1}$};
	\node at (1.5,0.5) {$\mathbf{1}$};
	\node at (2.5,0.5) {$\mathbf{1}$};
	\node at (3.5,0.5) {$\mathbf{1}$};
	\node at (4.5,0.5) {$\mathbf{1}$};
	\node at (5.5,0.5) {$\mathbf{1}$};
	\node at (6.5,0.5) {$\mathbf{0}$};
	\node at (7.5,0.5) {$\mathbf{0}$};
	\end{tikzpicture}
	\caption{Two skew overlaps between $Z$ and an element $A\in L_{m\times n}$}
	\label{AZ}
\end{figure}

\noindent \emph{Proof.} Let $A\in L_{m\times n}$. Vertical overlaps between $A$ and $Z$ do not create control windows with identical corresponding entries since, for $i=2,3,\ldots,m$, it is $A_1\neq Z_i$ and $Z_1\neq A_i$ according to Definitions \ref{L_m_n_even} (or \ref{L_m_n_odd}) and the construction of $Z$. In case of skew overlaps (see Figure \ref{AZ}),  we observe that a prefix or a suffix of $A_1=Z_1$ is always involved in the control window and it is superimposed to a suffix or a prefix of $X$ or $A_i$ ($i=1,2,\ldots,m$).

If the involved rows are $A_1=Z_1$ and $X$, no their prefix or suffix can match since $A_1$ and $X$ are non-overlapping strings. If the involved row are $Z_1$ and $A_i$, then the proof is similar to the one of Proposition \ref{n_even} (or \ref{n_odd}) recalling that $X$ ends with $0$.  

\cvd  

The construction of the string $X$ strictly depends on the first row $A_1$ of an element of $L_{m\times n}$. For example, if $A_1=11010100$, the string $X=11111100$ can be used to construct the matrix $Z$ in order to expand the set $L_{m\times n}$:

\begin{center}
	\begin{tikzpicture}[scale=0.75]
	\draw (0,0) rectangle (8,5);
	
	\node at (0.5,4.5) {$\mathbf{1}$};
	\node at (1.5,4.5) {$\mathbf{1}$};
	\node at (2.5,4.5) {$\mathbf{0}$};
	\node at (3.5,4.5) {$\mathbf{1}$};
	\node at (4.5,4.5) {$\mathbf{0}$};
	\node at (5.5,4.5) {$\mathbf{1}$};
	\node at (6.5,4.5) {$\mathbf{0}$};
	\node at (7.5,4.5) {$\mathbf{0}$};
	
	\node at (-1,2.5) {$Z=$};

	\draw (0.2,4) -- (7.8,4);

	\node at (0.5,3.5) {$\mathbf{1}$};
	\node at (1.5,3.5) {$\mathbf{1}$};
	\node at (2.5,3.5) {$\mathbf{1}$};
	\node at (3.5,3.5) {$\mathbf{1}$};
	\node at (4.5,3.5) {$\mathbf{1}$};
	\node at (5.5,3.5) {$\mathbf{1}$};
	\node at (6.5,3.5) {$\mathbf{0}$};
	\node at (7.5,3.5) {$\mathbf{0}$};
	
	\node at (0.5,2.5) {$\mathbf{1}$};
	\node at (1.5,2.5) {$\mathbf{1}$};
	\node at (2.5,2.5) {$\mathbf{1}$};
	\node at (3.5,2.5) {$\mathbf{1}$};
	\node at (4.5,2.5) {$\mathbf{1}$};
	\node at (5.5,2.5) {$\mathbf{1}$};
	\node at (6.5,2.5) {$\mathbf{0}$};
	\node at (7.5,2.5) {$\mathbf{0}$};
	
	\node at (0.5,1.5) {$\mathbf{1}$};
	\node at (1.5,1.5) {$\mathbf{1}$};
	\node at (2.5,1.5) {$\mathbf{1}$};
	\node at (3.5,1.5) {$\mathbf{1}$};
	\node at (4.5,1.5) {$\mathbf{1}$};
	\node at (5.5,1.5) {$\mathbf{1}$};
	\node at (6.5,1.5) {$\mathbf{0}$};
	\node at (7.5,1.5) {$\mathbf{0}$};

	\node at (0.5,0.5) {$\mathbf{1}$};
	\node at (1.5,0.5) {$\mathbf{1}$};
	\node at (2.5,0.5) {$\mathbf{1}$};
	\node at (3.5,0.5) {$\mathbf{1}$};
	\node at (4.5,0.5) {$\mathbf{1}$};
	\node at (5.5,0.5) {$\mathbf{1}$};
	\node at (6.5,0.5) {$\mathbf{0}$};
	\node at (7.5,0.5) {$\mathbf{0}$};

	\end{tikzpicture}
\end{center}

\noindent
Nevertheless, if $A_1=11100010$, then the same $X=11111100$
can not be used to construct $Z$ since its suffix $11100$ matches the prefix of length $5$ of $A_1$.

Unfortunately, once fixed a string $A_1$, we are not able to find a general criterion to construct a suitable string $X$.

\section{Further developments}
A first attempt to improve this study could start from the last hint of the preceding section. If $A_1=1u0$ with $u\in D_{n-2}$ is a \texttt{type} $\alpha$ Dyck word of length $n$, is there always a string $X$ with the required properties? In the positive case, is there a general procedure to construct it? If this procedure were found, it could shed light on the possibility of developing the set $L_{m\times n}$ in order to make it a non-expandable set.

\medskip
As far as the expandability is concerned, we notice that the types of rows of a matrix from the second to the last one depends on the first row $A_1$. For example,
if $A_1=1^{n/2}0^{n/2}$, then a row $A_i$ ($i\neq 1$) can be the following: $A_i=1^4v$ with $v\in D_{n-4}$ and $v\neq 1^{(n-4)/2}0^{(n-4)/2}$. This is a different row type with respect to \texttt{type 1, 2, 3, 4, 5}.  
It is not difficult to show that this new row and $A_1$ are non overlapping strings, so that it can be used to generate new matrices in $L_{m\times n}$.
Consequently, we could ask to find all the possible strings which are non-overlapping with $A_1$. Again, is there a general and efficient procedure which, fixed $A_1$, lists all the non-overlapping strings with it?

\medskip
Finally, an attempt aimed to the construction of variable-dimension matrices which are non-overlapping could be carried on, following up what already done \cite{bilosolo} for strings set.

\section{Acknowledgements}

This work has been partially supported by the GNCS project 2018 ``Propriet\`a combinatorie e rilevamento di pattern in 
strutture discrete lineari e bidimensionali''.


\end{document}